\documentclass[mathpazo]{cicp}

\usepackage{amsmath,epsf,cite}
\usepackage{amssymb,amsthm,tocvsec2}
\usepackage{bigints}
\usepackage{graphicx}
\usepackage{sidecap}
\usepackage{wrapfig}
\usepackage{pict2e}
\usepackage{float}
\usepackage{graphicx}
\usepackage{setspace}
\usepackage{float}
\usepackage{color}
\usepackage{mathrsfs}
\usepackage{epsfig,epsf,latexsym,subfigure}
\usepackage{fancybox}
\usepackage{fancyhdr}
\usepackage{mdwlist}
\usepackage{paralist}
\usepackage{setspace}
\usepackage{tikz}
\usetikzlibrary{shapes,arrows}
\usepackage{enumitem}
\setlist{nolistsep}
\usepackage{indentfirst}
\usepackage{epstopdf}

\numberwithin{equation}{section}
\numberwithin{figure}{section}
\numberwithin{table}{section}

\theoremstyle{plain}

\theoremstyle{remark}

\def\cG{\mathcal{G}}

\def\bx{\mathbf{x}}
\def\cN{\mathcal{N}}

\def\hR{\mathbb{R}}

\newcommand{\ben}{\begin{eqnarray}}
\newcommand{\een}{\end{eqnarray}}
\newcommand{\benx}{\begin{eqnarray*}}
\newcommand{\eenx}{\end{eqnarray*}}
\newcommand{\beq}{\begin{equation}}
\newcommand{\eeq}{\end{equation}}
\newcommand{\beqx}{\begin{equation*}}
\newcommand{\eeqx}{\end{equation*}}
\newcommand{\bea}{\begin{array}}
\newcommand{\eea}{\end{array}}
\newcommand{\bef}{\begin{figure}[H][H]}
\newcommand{\eef}{\end{figure}}
\newcommand{\be}{\begin{equation}}
\newcommand{\ee}{\end{equation}}

\newcommand{\bse}{\begin{subequations}}
\newcommand{\ese}{\end{subequations}}

\def\hR{\mathbb R}

\def\be{\mathbf{e}}

\def\bx{\mathbf{x}}

\def\bD{\mathbf{D}}

\begin{document}
\title{Discovering Phase Field Models from Image Data with the Pseudo-spectral Physics Informed Neural Networks}

\author[J. Zhao]{Jia Zhao \corrauth}
\emails{ {Department of Mathematics \& Statistics, Utah State University, Logan, UT, USA \tt jia.zhao@usu.edu.} (J.~Zhao)}

\begin{abstract}
In this paper, we introduce a new deep learning framework for discovering the phase field models from existing image data. The new framework embraces the approximation power of physics informed neural networks (PINN), and the computational efficiency of the pseudo-spectral methods, which we named pseudo-spectral PINN or SPINN. Unlike the baseline PINN, the pseudo-spectral PINN has several advantages. First of all, it requires less training data. A minimum of two snapshots with uniform spatial resolution would be adequate. Secondly, it is computationally efficient, as the pseudo-spectral method is used for spatial discretization. Thirdly, it requires less trainable parameters compared with the baseline PINN. Thus, it significantly simplifies the training process and assures less local minima or saddle points. We illustrate the effectiveness of pseudo-spectral PINN through several numerical examples. The newly proposed pseudo-spectral PINN is rather general, and it can be readily applied to discover other PDE-based models from image data.  
\end{abstract}


\maketitle

\section{Introduction}
The phase field models have been widely appreciated for investigating moving interface or multiphase problems. Among them, the well-known models are the Allen-Cahn equation (AC)  and the Cahn-Hilliard (CH) equation \cite{Cahn&H1958}. In general, most existing phase field models are thermodynamically consistent, i.e., they satisfy the first and second thermodynamic laws. In particular, for most scenarios, the temperature is assumed constant, such that the Helmholtz free energy is non-increasing in time.  Therefore, Denoting the state variable as $\phi$ and the free energy as $E$, the phase field models can usually be written in a generic gradient flow form
\beq \label{eq:gradient-flow}
\partial_t \phi = \cG \frac{\delta E}{\delta \phi},
\eeq 
with proper initial values and boundary conditions. Here $\cG$ is a semi negative-definite operator known as the mobility operator. In general, the free energy $E$ contains a conformation entropy and a bulk potential part. In this paper, we consider the free energy
\beq \label{eq:free-energy}
E = \int_\Omega \Big[ \frac{\varepsilon^2}{2} |\nabla \phi|^2 + F(\phi) \Big] dx.
\eeq 
with $\varepsilon$ a model parameter. We use $f(\phi) = F'(\phi)$ to denote the bulk chemical potential.
For the Allen-Cahn equation, the mobility is usually chosen as $\cG=-M$, with $M>0$ a constant. Then, the Allen-Cahn equation reads
\beq \label{eq:AC}
\partial_t \phi = M ( \varepsilon^2 \Delta \phi - f(\phi)).
\eeq 
Similarly, for the Cahn-Hilliard equation, the mobility is usually chosen as $\cG= M \Delta$, with $M>0$ a constant. Then the Cahn-Hilliard equation reads as
\beq \label{eq:CH}
\partial_t \phi = M \Delta ( -\varepsilon^2 \Delta \phi + f(\phi)).
\eeq 
Here the choice of the bulk potential $F(\phi)$ will depend on the material properties, and the double-well potential and Flory-Huggins free energy are the common choices.

Though there is a considerable amount of work on solving the phase field models \cite{Wang&Wang&WiseDCDS2010,Shen10_1, Han&WangJCP2014, Guillen&Tierra2014,Yang&Zhao&WangJCP2017,Gong-HIEQSAV}, there is less achievement on solving the inverse problem, i.e., discovering the phase field models from data.
In the classical approach,  the PDE models are usually derived based on empirical observation. The PDE models with free parameters are introduced first, and the free parameters will be fitted with data afterward. In this paper, we introduce a novel method using the physics informed neural network approach, where the model would be directly learned from data. In particular, we assume the bulk potential $F(\phi)$ is unknown, and we will discover it via learning from the existing data using the deep neural network.

The deep neural network has been widely used to investigate problems in various fields. Mathematically, the feed-forward neural network could be defined as compositions of nonlinear functions. Given an input $x \in \hR^{n_1}$, and denote the output of the $l$-th layer as $a^{[l]} \in \hR^{n_l}$, which is the input for $l+1$-th layer. In general, we can define the neural network as \cite{DeepLearningMath} 
\beq \label{eq:NN}
\bea{l}
a^{[1]} = x \in \hR^{n_1}, \\
a^{[k]} = \sigma \Big( W^{[k]} a^{[k-1]} + b^{[k]} \Big) \in \hR^{n_k}, \quad \mbox{ for } k=2,3,\cdots, L,
\eea 
\eeq 
where  $W^{[k]} \in \hR^{n_k \times n_{k-1}}$ and $b^{[k]} \in \hR^{n_k}$ denote the weights and biases at layer $k$ respectively, $\sigma$ denotes the activation function.
Many interesting work have been published in terms of solving or discovery of partial differential equations \cite{RaissiJCP,DeepXDE,DHPM,XiuWu1,Brunton3932,RungeKutta-NeuralNetwork,HanEPNAS2018,E,complex-geometries,Yu, PRL-image}. Here we briefly discuss some relevant work of PDE discovery with deep neural network or machine learning method. In \cite{RaissiJCP}, the authors assume the PDE model is a linear combination of a few known terms, and try to fit the coefficient from data. Based on the format of known data, they propose two approaches for calculating the loss function. In \cite{SINDY-PDE}, the authors assume the PDE model is a linear combination of terms in a predefined dictionary (with many terms), and use the sparse regularization to coefficients are zero. Thus, after fitting with the data, only a few non-zero terms will be left, and the model is discovered. In \cite{PDE-net}, the authors introduce a feed-forward deep neural network to predict the PDE dynamics and uncover the underlying PDE model (in a black box form) simultaneously.

However, in most of the existing work, there is strong assumption on the availability of data. For instance, most existing PDE discovery strategies request a well sampled data across time and spatial domain, and they usually assume the data collected at different time slots are from a single time sequence. In reality, especially from lab experiments, the data are usually sampled as snapshots (images) from various time slots, and the images might have limited spatial resolution, i.e. the solution values are only known at certain locations.  The major goal of this paper is to address the problem of PDE discovery with image snapshots. Therefore, we add to the efforts by introducing the pseudo-spectral PINN, which provides a handy approach to deal with data that are provided as image snapshots. And in particular, we focus on the phase field models alone.

The rest of of this paper is organized as follows. In Section 2, we introduce some notations, and set up the PDE discovery problem. In Section 3, we introduce the pseudo-spectral physics informed neural networks by introducing the network stricture and the various definitions of the loss function. Afterward, several examples are shown to demonstrate the approximation power of the pseudo-spectral PINN in Section 4. A brief conclusion is drawn in the last section.

\section{Problem Setup}
Recall from \eqref{eq:gradient-flow}, we use $\phi(\bx, t)$ to denote the state variables, where $\bx \in \Omega \subset \mathbb{R}^d$ are the spatial coordinates, and $t\in(0, T]$ denotes the time. The baseline PINN \cite{RaissiJCP} would require a rather amount of data pairs $\{(\bx_i, t_i, \phi(\bx_i, t_i))\}_{i=1}^N$ properly sampled from the domain $\Omega \times [0, T]$. Also, it requires all the data are from a single time sequence, given a deep neural network $\mathcal{N}(\bx, t)$ is needed to approximate $\phi(\bx, t)$ as an aiding neural network for the discovery the PDE.
However, this is usually not the case. In practice, the data are usually sampled as image snapshots at various time slots. For instance, for certain experiments, a picture is taken at $t_0$, and after $\delta t$, another picture is taken at $t_1=t_0+\delta t$.

To clearly describe the problem, we need to introduce some notations.
Let $N_x,N_y$ be two positive even integers. The spatial domain $\Omega = [0,L_x]\times[0,L_y]$ is uniformly partitioned with mesh size $h_x = L_{x}/N_{x}, h_y = L_{y}/N_{y}$ and
$$
\Omega_{h} =
\left\{(x_{j},y_{k})|x_{j} = j h_x, y_{k} = kh_y,~0\leq j\leq
N_{x}-1,0\leq k\leq N_{y}-1\right\}.
$$
In order to derive the algorithm conveniently, we denote the discrete gradient operator and the discrete Laplace operator 
$$\nabla_h = \left(
\begin{array}{c}
\bD_1^x\textcircled{x}\\
\bD_1^y\textcircled{y}
\end{array}\right),
\quad
\Delta_h = \nabla_h\cdot\nabla_h = (\bD_1^x)^2\textcircled{x}+(\bD_1^y)^2\textcircled{y},
$$
following our previous work \cite{Chen&Zhao&GongCICP2019}.

Let $V_h = \{ u|u ={u_{j,k} (x_j, y_k) \in \Omega_h} \}$ be the space of grid functions on $\Omega_h$.
For such cases, the data are collected in the form as
\beq \label{eq:general_data}
\{ (\Phi_i^{(1)}, \Phi_i^{(2)}, \delta_i) \}_{i=1}^N  \subset V_h \times V_h \times \hR^{+},
\eeq 
where $\Phi_i^{(1)} \in V_h$ and $\Phi_i^{(2)} \in V_h$ are two snapshots with a time lag $\delta_i>0$ between the two states, and $N$ is the total number of snapshot pairs. And the goal is to discover the bulk function $f$ in the phase field models, with the data collected in the form of \eqref{eq:general_data}.

We emphasis that $\delta_i$ is not required to be constant, removing the assumption in \cite{XiuDongbin}. Also, different data pairs $(\Phi_i^{(1)}, \Phi_i^{(2)})$ and $(\Phi_j^{(1)}, \Phi_j^{(2)})$ $i\neq j$ are not necessarily from the same single time sequence, removing the requirement of baseline PINN \cite{RaissiJCP}.  

\section{Pseudo-spectral Physics Informed Neural Networks}
\subsection{Neural Network Structure}

For simplicity of notations, we introduce our idea by applying it on a generic example.  Given the PDE problem
\beq \label{eq:generic-model}
\partial_t \phi = \cG \Big[ g(\phi, \nabla \phi, \Delta \phi) + f(\phi) \Big], \quad \bx \in \Omega, \quad t \in [0, T],
\eeq 
with periodic boundary condition, where operator $g$ is known, but operator $f$ is unknown. If we are provided with data pairs as \eqref{eq:general_data}, the major goal of this paper is to identify the operator/functional $f$, thus discover the PDE model in \eqref{eq:generic-model}. 

The generic way to solve it is by introducing a deep neural network 
\beq
\cN_f: \phi \rightarrow \cN_f(\phi; \theta),
\eeq 
to approximate $f(\phi)$, where $\theta$ represents the free parameters.  In this paper, we assume $f$ is only a function of $\phi$ for simplicity. Notice the idea introduced applies easily to cases where $f$ is a function of $\phi$ and its derivatives, saying $f:=f(\phi, \nabla \phi, \Delta \phi)$.

Since the data pairs \eqref{eq:general_data} are uniformly sampled and the periodic boundary condition is considered for the PDE problem \eqref{eq:generic-model}, it is advisable to simply the problem \eqref{eq:generic-model} by discretizing in space using pseudo-spectral method. Afterwards, the PDE problem in \eqref{eq:generic-model} can be written as
\beq \label{eq:generic-model-semi-discrete}
\partial_t \Phi_{ij} = \cG_h \Big[ g(\Phi_{ij}, \nabla_h \Phi_{ij}, \Delta_h \Phi_{ij}) + f(\Phi_{ij}) \Big], \quad i,j=1,2,\cdots,M, \quad t \in [0, T],
\eeq 
where $\cG_h$ is the spatially discretized mobility operator. Here we use $\Phi \in V_h$ to denote the discrete function values of $\phi$ on $\Omega_h$.
We emphasis the periodic boundary condition is discussed in this paper. If other type of boundary conditions is considered, finite-difference or finite-element method for spatial discretization might be more proper.

\begin{remark}
Notice, in the baseline PINN, an aiding neural network  $\cN: (\bx, t) \rightarrow \cN(\bx, t)$ is introduced to approximate $\Phi(\bx,t)$, which is computationally expensive and is applicable only to a single time sequence. It will be apparent that applying the physics informed neural networks on the semi-discrete problem \eqref{eq:generic-model-semi-discrete} will have several advantages.
\end{remark}

There are mainly two components (ingredients) in the deep neural network method. First of all, we shall define what the neural network meant to approximate (along with its structures, input, output, and hidden layers); define the loss function (which is enforced with known physics). Hence, for the problem above, we introduce a neural network $\mathcal{N}_f:\phi \rightarrow \mathcal{N}_f(\phi; \theta)$ to approximate $f$. Next, we need to define the loss function.


\subsection{Loss function}
We split $g$ in \eqref{eq:generic-model-semi-discrete} as  
$$
g(\Phi, \nabla_h \Phi, \Delta_h \Phi) = L_g(\Phi) + N_g(\Phi),
$$
where $L_g$ is the linear operator and $N_g$ is the rest, i.e.
\beq
\partial_t \Phi = \cG_h \Big[ L_g(\Phi) + N_g(\Phi) + f(\Phi) \Big].
\eeq 
To solve it numerically in the time interval $[t_i, t_i +\delta_i]$, we can propose the stabilized scheme
\beq
\frac{1}{\delta_i} ( \Phi_{t_i+\delta_i} - \Phi_{t_i} ) =\cG_h \Big[ L_g (\Phi_{t_i+\delta_i}) + N_g(\Phi_{t_i}) + f(\Phi_{t_i}) + C(\Phi_{t_i +\delta_i} - \Phi_{t_i}) \Big],
\eeq 
where $C$ is a stabilizing operator, which could be chosen as
\beq
C  = \sum_{i=0}^2 (-1)^i c_i (\Delta_h)^{2i},
\eeq 
with $c_i$'s are constants. Then, we have
\beq \label{eq:BDF1-scheme}
\Phi_{t_i+\delta_i} = (1 - \delta_i (\cG_h(C+L_g)))^{-1}\Big(\Phi_{t_i} + \delta_i \cG_h\Big[N_g(\Phi_{t_i}) + f(\Phi_{t_i}) - C \Phi_{t_i} \Big] \Big).
\eeq 
Notice the scheme \eqref{eq:BDF1-scheme} is first-order accurate in time. When $\delta_i$ is small enough, the expression for $\Phi_{t_i+\delta_i}$ is accurate.
Inspired by this, we can introduce our linear SPINN, to discover \eqref{eq:generic-model} from the data \eqref{eq:general_data}.

\begin{definition}[Linear SPINN loss function]
Denote the neural network $\cN_f: \Phi \rightarrow \cN_f(\Phi; \theta)$ where $\theta$ are the free parameters. Given $(\Phi_i^{(1)},\delta_i)$, we can approximate $\Phi_i^{(2)}$ via $\cN_R(\Phi_i^{(1)}, \delta_i)$, which is defined by
\beq
\cN_R: (\Phi_i^{(1)}, \delta_i) \rightarrow  (1 - \delta_i (\cG_h(C+L_g)))^{-1}\Big(\Phi_i^{(1)} + \delta_i \cG_h \Big[N_g(\Phi_i^{(1)}) + \cN_f(\Phi_i^{(1)};\theta) - C \Phi_i^{(1)} \Big] \Big).
\eeq 
Hence, the loss function is defined as
\beq
L(\theta) = \sum_{i=1}^N \| \Phi_i^{(2)} - \cN_R(\Phi_i^{(1)}, \delta_i; \theta) \|_2^2.
\eeq 
\end{definition}

When the time step $\delta_i$ is small, by minimizing $L(\theta)$, we can identify $\cN_f$, which in turn discover the PDE problem in \eqref{eq:generic-model}. 
For certain cases, $\delta_i$ might be large. A single step marching scheme might be inaccurate. We thus introduce a loss function that is based on a more general recursive Linear SPINN.

\begin{definition}[Recursive Linear SPINN loss function]
Denote the neural network $\cN_f: \Phi \rightarrow \cN_f(\Phi; \theta)$. Given $(\Phi_i^{(1)},\delta_i)$, we can approximate $\Phi_i^{(2)}$ via a mapping $\cN_{R_K}(\Phi_i^{(1)}, \delta_i)=R_K$,
where $R_K$ is defined recursively as
\beq
\bea{l}
R_0 = \Phi_i^{(1)}, \\
R_j = (1 - \delta_i (\cG_h(C+L_g)))^{-1}\Big(R_{j-1} + \delta_i \cG_h\Big[N_g(R_{j-1}) + f(R_{j-1}) - C R_{j-1} \Big] \Big),
j=1,2,\cdots,K,
\eea
\eeq 
with $K$ and $c_i$'s hyper-parameters. And the loss function is traced as
\beq
L(\theta) = \sum_{i} \| \Phi_i^{(2)} - \cN_{R_K}(\Phi_i^{(1)}, \delta_i; \theta) \|_2^2.
\eeq 
\end{definition}

\begin{remark}
Similarly, we can design the Neural network inspired by a second-order or higher-order numerical scheme. For instance, consider the following predictor-corrector second-order scheme.
To solve it numerically in the time interval $[t_i, t_i +\delta_i]$, we introduce the stabilized second-order scheme in two step: 
\begin{itemize}
\item 
First of all, we can obtain $\hat{\Phi}_{t_i+\frac{\delta_i}{2}}$ via
\beq
\frac{ \hat{\Phi}_{t_i+\frac{\delta_i}{2}} - \Phi_{t_i}}{\delta_i/2} = \cG_h \Big[ L_g (\hat{\Phi}_{t_i+\frac{\delta_i}{2}}) + N_g(\Phi_{t_i}) + f(\Phi_{t_i}) + C(\hat{\Phi}_{t_i +\frac{\delta_i}{2}} - \Phi_{t_i}) \Big]. \\
\eeq 
\item Next, we can obtain $\Phi_{t_i+\delta_i}$ via solving 
\beq
\frac{ \Phi_{t_i+\delta_i} - \Phi_{t_i}}{\delta_i} = \cG_h \Big[ L_g (\frac{\Phi_{t_i+\delta_i}+\Phi_{t_i}}{2}) + N_g(\hat{\Phi}_{t_i+\frac{\delta_i}{2}}) + f(\hat{\Phi}_{t_i+\frac{\delta_i}{2}}) + C(\frac{\Phi_{t_i +\delta_i}+\Phi_{t_i}}{2} - \hat{\Phi}_{t_i+\frac{\delta_i}{2}}) \Big],
\eeq 
\end{itemize}

Therefore, if we definite the neural network $\cN_f: \phi \rightarrow \cN_f(\phi; \theta)$, once given $(\Phi_i^{(1)}, \delta_i)$, we can approximate $\Phi_i^{(2)}$ via the mapping defined as
\beq
\bea{l}
\cN_R: (\Phi_i^{(1)}, \delta_i) \rightarrow (1 - \frac{\delta_i}{2} (\cG_h(C+L_g)))^{-1}\Big(\Phi_i^{(1)} \\
\qquad \qquad + \delta_i \cG_h \Big[L_g(\frac{\Phi_i^{(1)}}{2})  + N_g(\hat{\Phi}_{i+\frac{\delta_i}{2}}) + \cN_f(\hat{\Phi}_{i+\frac{\delta_i}{2}};\theta) + C(\frac{ \Phi_i^{(1)}}{2} -\hat{\Phi}_{i+\frac{\delta_i}{2}}) \Big] \Big),
\eea
\eeq 
where $\hat{\Phi}_{i+\delta_i/2}$ is defined by
\beq
\hat{\Phi}_{i+\frac{\delta_i}{2}} =  
(1 - \frac{\delta_i}{2} (\cG_h(C+L_g)))^{-1}\Big(\Phi_i^{(1)} + \frac{\delta_i}{2} \cG_h \Big[N_g(\Phi_i^{(1)}) + \cN_f(\Phi_i^{(1)};\theta) - C \Phi_i^{(1)} \Big] \Big).
\eeq  
Then the loss function could be defined similarly. Some other ideas, such as exponential time integration can also be utilized for designing neural network structure. For brevity, these ideas will not be pursued in this paper. Interested readers are encouraged to explore these interesting topics.
\end{remark}

As an analogy, we can mimic the Runge-Kutta method for solving \eqref{eq:generic-model} to design neural networks.Inspired by the idea of four-stage explicit Runge-Kutta method for the time discretization, we can introduce the following loss function.
\begin{definition}[Four-stage explicit Runge-Kutta SPINN loss function]
Denote the neural network $\cN_f: \Phi \rightarrow \cN_f(\Phi; \theta)$. Given $(\Phi_i^{(1)}, \delta_i)$,  we can approximate $\Phi_i^{(2)}$ by the mapping $\cN_R(\Phi_i^{(1)}, \delta_i)$ defined by the four stage method
\beq
\cN_R: (\Phi_i^{(1)}, \delta_i) \rightarrow \Phi_i^{(1)} + \frac{\delta_i}{6} (K_1+2K_2+2K_3+K_4),
\eeq 
where
\beq
\bea{l}
K_1 =  \cG_h \Big[ g(\Phi_i^{(1)}) + \cN_f(\Phi_i^{(1)}; \theta) \Big], \\
K_2 =  \cG_h \Big[ g(\Phi_i^{(1)} + \frac{\delta_i}{2} K_1) + \cN_f(\Phi_i^{(1)} + \frac{\delta_i}{2} K_1; \theta)\Big],  \\
K_3 =  \cG_h \Big[g(\Phi_i^{(1)}+\frac{\delta_i}{2} K_2) + \cN_f(\Phi_i^{(1)}+\frac{\delta_i}{2} K_2; \theta)\Big],  \\
K_4=  \cG_h \Big[g(\Phi_i^{(1)} +\delta_i K_3) + N_f(\Phi_i^{(1)} +\delta_i K_3;\theta)\Big].
\eea 
\eeq 
Then the loss function can be defined as
\beq
L(\theta) = \sum_{i=1}^N \| \Phi_i^{(2)} - \cN_R(\Phi_i^{(1)}, \delta_i)\|^2. 
\eeq 
\end{definition}

Similarly,  when the time step $\delta_i$ is large, we can define the loss function via the explicit RK SPINN recursively as below. 
\begin{definition}[Recursive four-stage Runge-Kutta SPINN loss function]
Denote the neural network $\cN_f: \Phi \rightarrow \cN_f(\Phi; \theta)$. Given the data $(\Phi_i^{(1)}, \delta_i)$, we can approximate $\Phi_i^{(2)}$ via the $K$-step recursive mapping denoted as $\cN_{R_K}: (\Phi_i^{(1)}, \delta_i) \rightarrow R_K$, where
\beq
\bea{l}
R_0 = \Phi_i^{(1)}, \\
R_j = R_{j-1} + \frac{\delta_i}{6K} (K_1^{j-1}+2K_2^{j-1}+2K_3^{j-1}+K_4^{j-1}), j=1,2,\cdots, K,  \\
K_1^{j-1} = \cG_h\Big[ g(R_{j-1}) + \cN_f(R_{j-1}; \theta) \Big],   \\
K_2^{j-1} =  \cG_h\Big[ g(R_{j-1} + \frac{\delta_i}{2K} K^{j-1}_1) + \cN_f(R_{j-1}+ \frac{\delta_i}{2K} K^{j-1}_1; \theta)\Big],  \\
K_3^{j-1} =  \cG_h\Big[g(R_{j-1}+\frac{\delta_i}{2K} K^{j-1}_2) + \cN_f(R_{j-1}+\frac{\delta_i}{2K} K^{j-1}_2; \theta)\Big],  \\
K_4^{j-1}= \cG_h\Big[ g(R_{j-1} +\frac{\delta_i}{K} K^{j-1}_3) + N_f(R_{j-1} +\frac{\delta_i}{K} K^{j-1}_3;\theta)\Big].
\eea 
\eeq 
And the loss function  is  defined by
\beq
L(\theta) = \sum_{i=1}^N \| \Phi_i^{(1)} - \cN_{R_K}(\Phi_i^{(1)}, \delta_i)\|^2 . 
\eeq 
\end{definition}

\begin{remark}
When the time step is even larger, one can use the implicit Runge-Kutta method for time discretization. However, given the intermediate stages are unknown, we need to introduce an extra neural network to approximate it:
\beq
\cN_q: (\Phi_i^{(1)}, \delta_i) \rightarrow  (\cN_{c_1}(\Phi_i^{(1)}, \delta_i), \cN_{c_2}(\Phi_i^{(1)}, \delta_i), \cdots, \cN_{c_q}(\Phi_i^{(1)}, \delta_i)),
\eeq 
where $\cN_{c_j}(\Phi_i^{(1)}, \delta_i)$ is to approximate $\Phi(t_i+c_j\delta_i)$. By following the idea in \cite{RaissiJCP}, we can define the following loss function 
\beq
\bea{l}
L(\theta) = 
\sum_{j=1}^q \Big\| \Phi_i^{(1)} -  \Big[ \cN_{c_j}(\Phi_i^{(1)}, \delta_i) + \delta_i \sum_{k=1}^q a_{jk} \cN_{c_k}(\Phi_i^{(1)}, \delta_i) \Big] \Big\|^2 \\
+ \sum_{j=1}^q \Big\| \Phi_i^{(2)} - \Big[ \cN_{c_j}(\Phi_i^{(1)}, \delta_i) + \delta_i \sum_{k=1}^{q} (a_{jk}-b_k) \cN_{c_k}(\Phi_i^{(1)}, \delta_i) \Big] \Big\|^2.
\eea 
\eeq 
For such a case, it requires expensive computational costs, and we will not investigate it in the current paper.
\end{remark}

\section{Numerical examples}


Next, we investigate the proposed SPINN with several examples. To identify the phase field models from data snapshots.

Recall that we assume $f(\phi)$ is unknown, and the goal is to identify it via learning the existing data in the form of \eqref{eq:general_data}. As explained, we define a neural network.
$
\cN_f:(\phi) \rightarrow \cN_f(\phi; \theta)
$ 
to approximate $f$ and use the loss functions as defined in the previous section. In the rest of this paper, we assume $\cN_f$ a feed-forward neural network, with $2$-hidden layers, with each hidden layer having $20$ neurons. The $\tanh$ activation function is applied in both hidden layers. During the training process, the $Adam$ method with default learning rate is used for 10,000 training iterations, followed by an L-BFGS-B optimization training process. The algorithms are implemented with Tensorflow.

For simplicity of discussion, in the all the examples below, $rand()$ generates random numbers between $[-1, 1]$ we chosen the domain $\Omega=[-1, 1]^2$, and choose $N_x=N_y=128$ in $\Omega_h$, i.e. the collected data $\Phi_i^{(1)}, \Phi_i^{(2)}$ are matrices in $\hR^{128,128}$. And we solve the PDE first with high-order-accurate scheme with uniform time steps and uniform spatial discretization. The numerical solutions are randomly sampled at different time slots as training data to inversely discover the bulk function $f$.  Given the free energy in \eqref{eq:free-energy}, we get $g(\phi) = -\varepsilon^2 \Delta \phi$, and $L_g(\Phi)=-\varepsilon^2 \Delta_h \Phi$. And we will chose $C=2$ for AC equation, and $C=-2\Delta_h$ for CH equation.

%

\begin{figure}[H]
\centering
\includegraphics[width=0.24\textwidth]{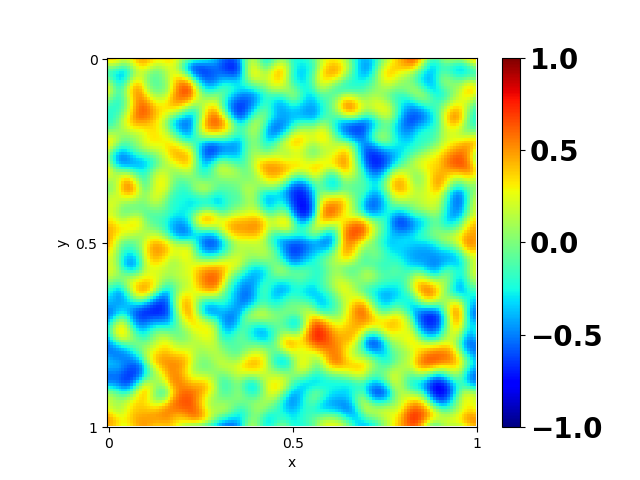}
\includegraphics[width=0.24\textwidth]{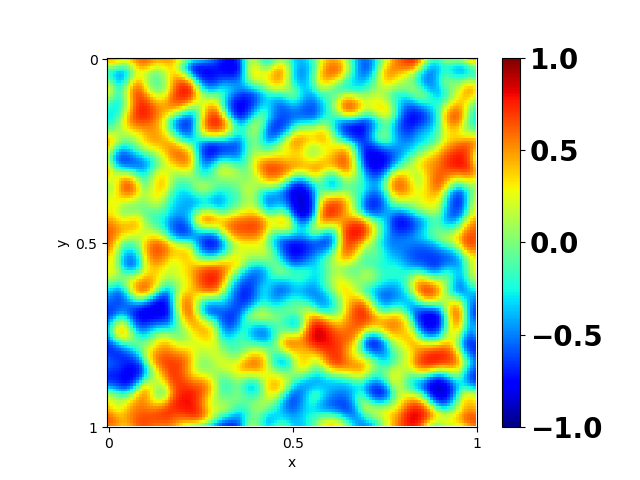}
\includegraphics[width=0.24\textwidth]{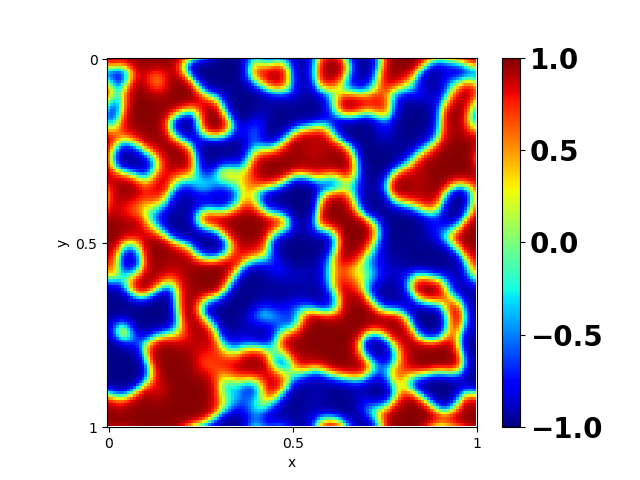}
\includegraphics[width=0.24\textwidth]{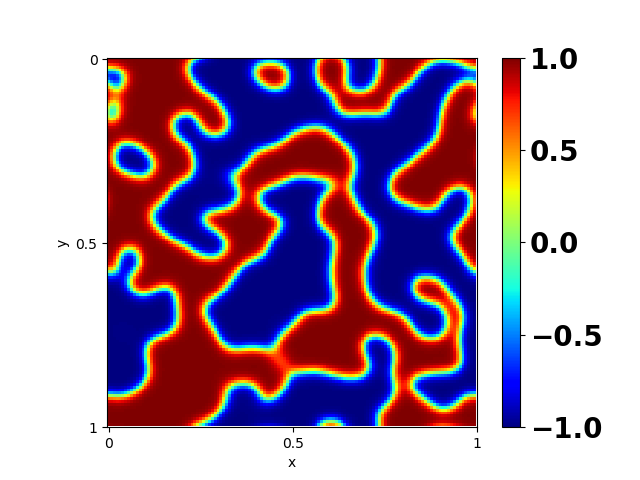}
\caption{Snapshots of the solution $\Phi$ for the Allen-Cahn equation at various time slots $t=0.325, 0.375, 0.575, 0.825$.}
\label{fig:AC2d-snapshots}
\end{figure}

\textbf{Example 1}.
In the first example, we generate data by solving the Allen-Cahn equation in \eqref{eq:AC} with $F(\phi) = \frac{1}{4}(1-\phi^2)^2$, which means $f(\phi) = \phi^3 - \phi$. And the parameters used are
$\varepsilon= 0.02$, $M=10$, with the initial condition $\phi(\bx, t=0) = 0.25 rand(\bx)$. Some snapshots of $\phi$ at various time slots are shown in Figure \ref{fig:AC2d-snapshots}.

\begin{figure}[H]
\centering
\includegraphics[width=0.9\textwidth]{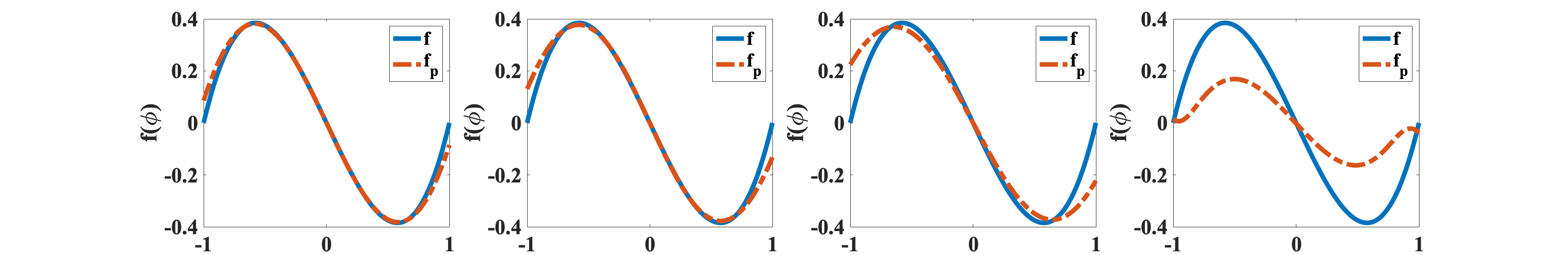}
\caption{Bulk function $f$ predicted using the Runge-Kutta SPINN loss function with $\delta =0.05, 0.1, 0.25, 0.5$.}
\label{fig:RKPINN-delta}
\end{figure}

First of all, we test out the Runge-Kutta SPINN approach. We choose $N=1$, i.e. use only a single data pair $(\Phi^{(1)}, \Phi^{(2)}, \delta)$ to train the neural network. We randomly choose $\Phi^{(1)}$, and study how the size of $\delta$ would affect the learned result. One experiment result with $\Phi^{(1)}$ chosen at $t=0.325$  and $\delta =0.05, 0.1, 0.25, 0.5$ are summarized in Figure \ref{fig:RKPINN-delta}. We observe that when the two snapshots are close, i.e. $\delta$ is small, the Runge-Kutta SPINN approach can accurately  learn the bulk function $f$. However, when the time step $\delta$ is large, its accuracy drops.

\begin{figure}[H]
\centering
\includegraphics[width=0.98\textwidth]{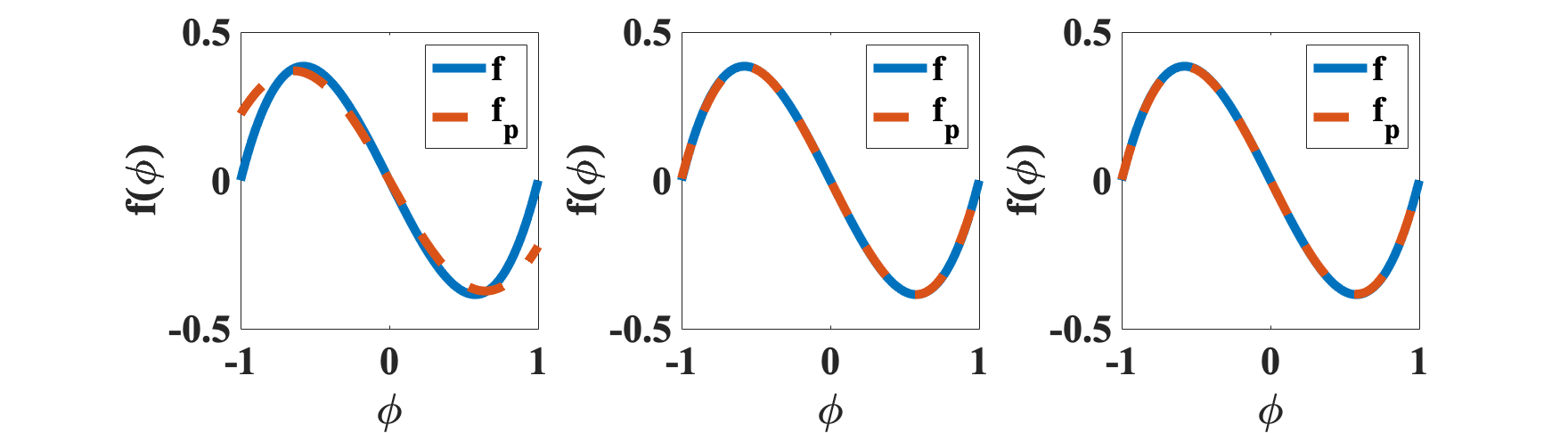}
\caption{Bulk function $f$ predicted using the recursive Runge-Kutta SPINN loss function, where $K=1,50,100$. This figure indicates with sufficient recursive stages, it gives accurate prediction of $f$. }
\label{fig:RKPINN-K}
\end{figure}

To overcome the inaccuracy when $\delta$ is large, we utilize the recursive Runge-Kutta SPINN approach. Here we fix the time step $\delta =0.25$, and test the effect by using different recursive stage $K=1,20,50,100$. The results are summarized in Figure \ref{fig:RKPINN-K}. We observe that, by increasing the recursive stages, the accuracy improves, and the bulk function $f$ can be learned accurately when sufficient stage is used.

\begin{figure}[H]
\centering
\includegraphics[width=0.95\textwidth]{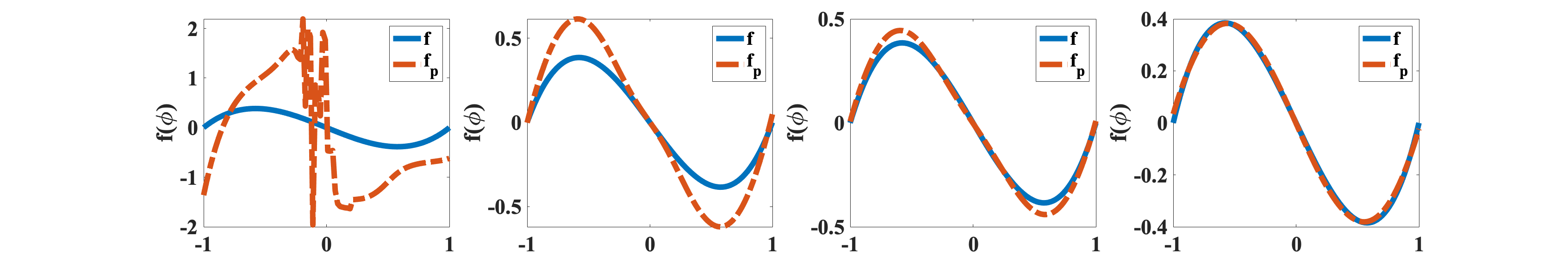}
\caption{Predicted bulk function $f$ with K=1,10,50,100. It shows ,with $K$ increasing, the accuracy of the SPINN improves. }
\label{fig:IPINN-K}
\end{figure}

\textbf{Example 2}.
However, when the time step is large enough, the recursive Runge-Kutta SPINN approach will not provide accurate approximate to $f$. For instance, with the time step $\delta =0.5$, the recursive Runge-Kutta SPINN approach fails. Meanwhile, the linear SPINN approach shows superior accuracy. As an example, we use a single data pair$(\Phi^{(1)}, \Phi^{(2)}, \delta)$  with a fixed time step $\delta =0.5$. We vary the recursive stage $K$ for the recursive linear SPINN, and the results of the learned function $f$ are summarized in Figure \ref{fig:IPINN-K}.
We observe that, even when the time step $\delta$ is large, the recursive linear SPINN approach still provides accurate approximate to $f$, so long as the recursive stage $K$ is large enough. 
We remark that the training strategy and the quality of training data might also be factors for the approximation accuracy, which we will not pursue in detail. Interested readers are strongly encouraged to explore.

\begin{figure}[H]
\centering
\subfigure[snapshots of $\phi$ at $t=0.325, 0.375$]{
\includegraphics[width=0.24\textwidth]{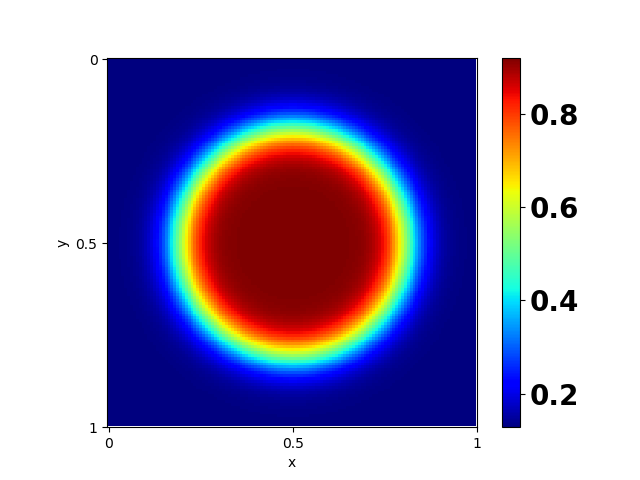}
\includegraphics[width=0.24\textwidth]{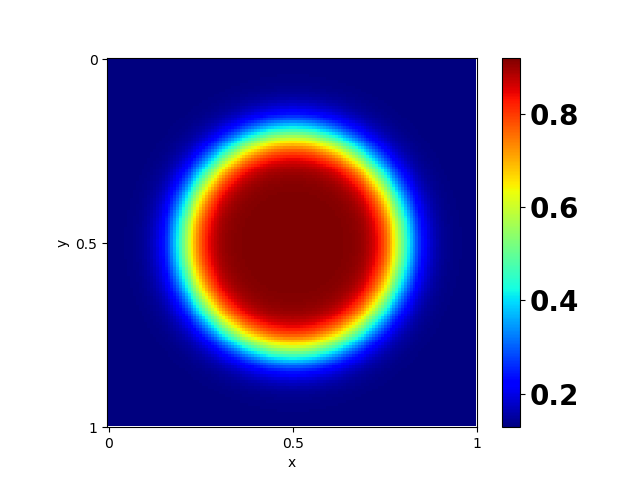}}
\subfigure[predicted $f$]{\includegraphics[width=0.24\textwidth]{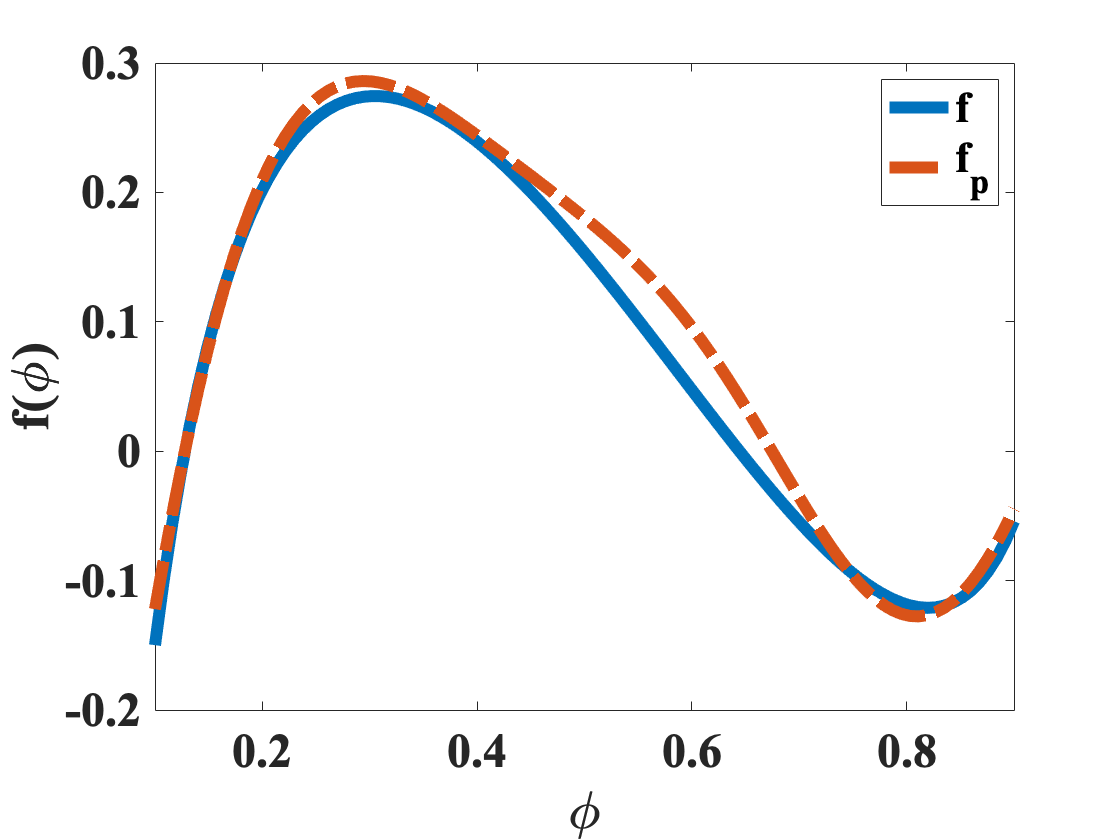}}
\caption{Predicted bulk function $f$ for the AC equation with Flory-Huggins bulk free energy using the Linear SPINN with $K=10$. (a) snapshots of the solution $\phi$ for the Allen-Cahn equation at time $t=0.325, 0.375$. (b) the predicted bulk function $f$ learned from the data in (a).}
\label{fig:AC2d-FH}
\end{figure}

\textbf{Example 3}. In the next example, we increase the problem complexity to identify a highly nonlinear bulk function. In details, we get the data by solving an Allen-Cahn equation with the Flory-Huggins free energy 
$F(\phi) = \phi \ln(\phi) + 0.5(1-\phi)\ln (1-\phi) + 2 \phi(1-\phi)$,
which means the bulk function $f = \ln \phi - 0.5 \ln(1-\phi) +1.5-2\phi$. The parameters used are $M=10$, $\varepsilon=0.1$, along with the initial condition $\phi(\bx,t=0)= \frac{1}{2} ( 1 + \tanh \frac{0.8 - \sqrt{x^2+y^2}}{\sqrt{2}\varepsilon})$. We randomly sample two snapshots with $\delta t = 0.05$, and train the neural network. A example of using two snapshots at $t=0.325, 0.375$ are shown in Figure \ref{fig:AC2d-FH}, where the two snapshots are shown in Figure \ref{fig:AC2d-FH}(a), and the predicted function $f$ is shown in Figure \ref{fig:AC2d-FH}(b). 
We observe the linear SPINN approach can learn the bulk function $f$ from only two images accurately.

\begin{figure}[H]
\centering
\includegraphics[width=0.24\textwidth]{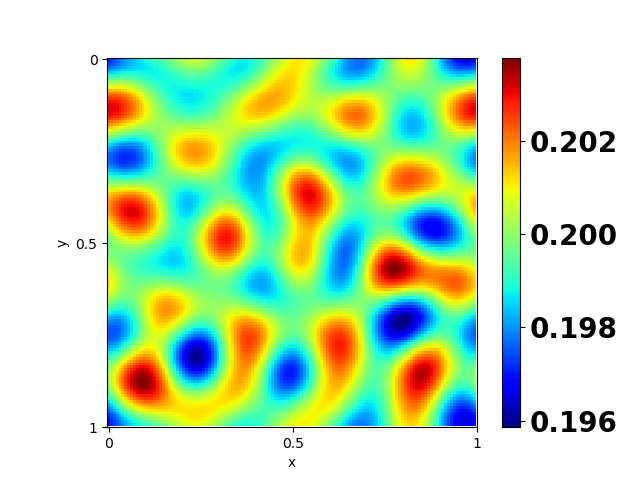}
\includegraphics[width=0.24\textwidth]{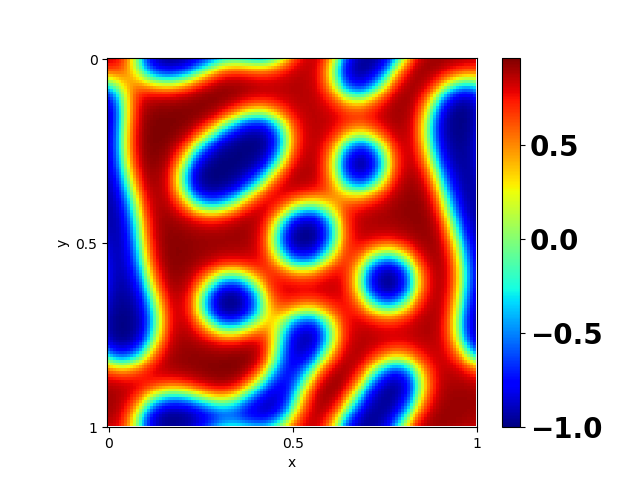}
\includegraphics[width=0.24\textwidth]{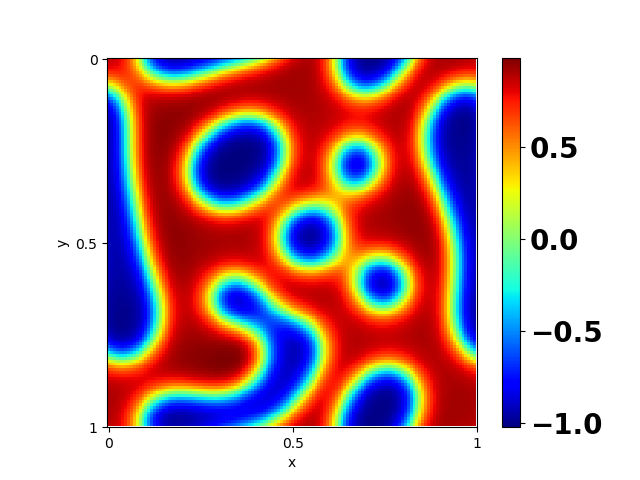}
\includegraphics[width=0.24\textwidth]{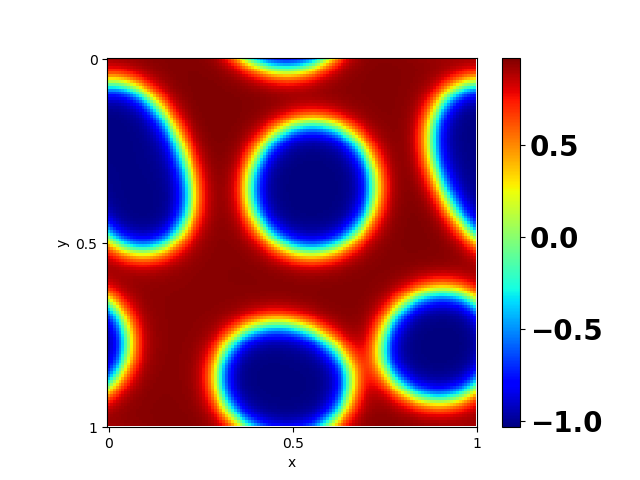}
\caption{Snapshots of the solution $\Phi$ for the Cahn-Hilliard equation at various time slots $t=0.05, 0.325, 0.375, 1$.}
\label{fig:CH-snapshots}
\end{figure}

\begin{figure}[H]
\centering
\includegraphics[width=0.9\textwidth]{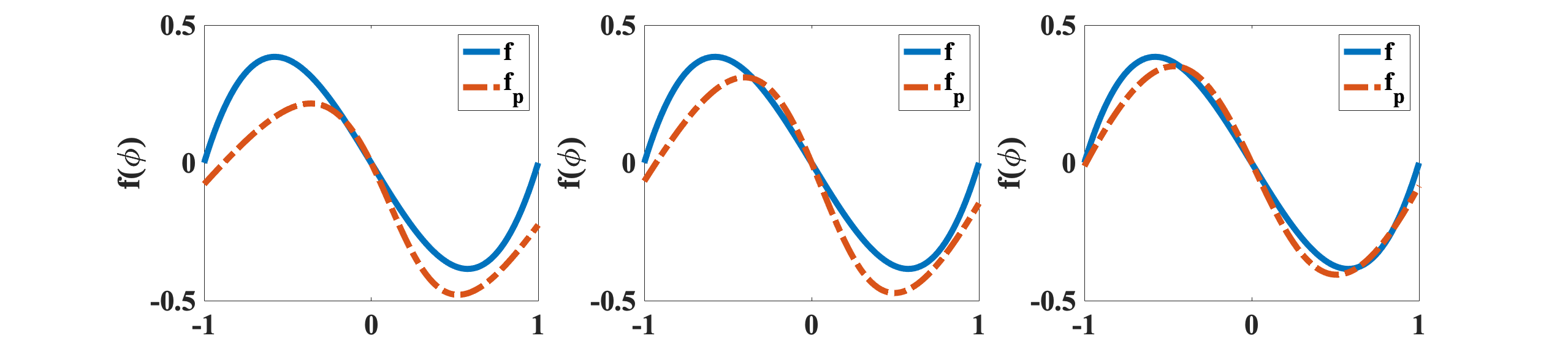}
\caption{Predicted bulk function $f$ for the Cahn-Hilliard equation with various data pairs $N=1,5, 10$.}
\label{fig:CH-f}
\end{figure}

\textbf{Example 4}.
In the last example, we use the linear SPINN approach to discover the bulk function $f$ from the solution snapshots of the Cahn-Hilliard equation in \eqref{eq:CH}. Here the data is obtained by solving the Cahn-Hilliard equation with $F(\phi) = \frac{1}{4} (\phi^2  -1)^2$, i.e. $f(\phi) = \phi^3 - \phi$, $M =1$ $\varepsilon=0.05$, and initial condition $\phi(t=0)=0.2+0.001 rand(\bx)$. For the loss function, we add an extra term $10^3\| \cN_f(0) \|^2$, to enforce the uniqueness of $\cN_f$. 
Some temporal snapshots for $\Phi$ are summarized in Figure \ref{fig:CH-snapshots}.
We randomly choose data $\{ (\Phi_i^{(1)}, \Phi_i^{(2)}, \delta_i) \}_{i=1}^N$, with $N=1, 5, 10$ data points with fixed time step $\delta_i = 0.05$. The learned result $f$ is summarized in Figure \ref{fig:CH-f}. We obverse the predicted bulk function $f$ has improved accuracy with more data used to train the neural network.

\section{Conclusion}
In this paper, we introduce a pseudo-spectral physics informed neural network to discover the bulk function in the phase field models. This newly proposed method well fits the common data collection strategy, i.e., taking snapshots/images at various time slots. The effectiveness of the proposed pseudo-spectral PINN, or SPINN, has been verified through identifying the bulk function $f$  of several phase field models. The idea of pseudo-spectral PINN is rather general, and it can be readily applied to discover other PDE models from image data, which will be investigated in our later research projects.

\section*{Acknowledgment}
Jia Zhao would like to acknowledge the support from NSF DMS-1816783.

\section*{Conflict of Interest Statement}
On behalf of all authors, the corresponding author states that there is no conflict of interest.

\bibliographystyle{unsrt}

\end{document}